\title{A note on graphs without short even cycles}

\author{
Thomas Lam
\thanks{Department of Mathematics, Massachusetts Institute of Technology,
77 Massachusetts Ave., Cambridge, MA 02139, USA.
E-mail: thomasl@math.mit.edu}
\and
Jacques Verstra\"{e}te
\thanks{Faculty of Mathematics, University of Waterloo,
200 University Avenue West, Waterloo, Ontario, Canada.
E-mail: jverstraete@math.uwaterloo.ca}
}

\documentclass [12pt]{article}
\usepackage{amsfonts}
\usepackage{amsmath}
\usepackage{graphicx}
\input labelfig.tex
\newtheorem{theorem}{Theorem}
\newtheorem{proposition}[theorem]{Proposition}

\newtheorem{lemma}[theorem]{Lemma}

\newenvironment{proof}{{\bf Proof.}}{\hfill{ }\vrule height10pt
width5pt depth1pt \medskip}

\newcommand{\floor}[1]{\lfloor #1 \rfloor }

\newcommand{\set}[1]{ \left\{ #1 \right\} }
\newcommand{\W}{\mathcal{W}}
\newcommand{\p}{\mathcal{P}}
\newcommand{\q}{\mathcal{Q}}
\newcommand{\brac}[1]{\left( #1 \right)}

\begin{document}

\maketitle

\vspace{-0.2in}

\begin{abstract}
In this note, we show that any $n$-vertex graph
without even cycles of length at most $2k$ has at most
$\frac{1}{2}n^{1 + 1/k} + O(n)$ edges, and
polarity graphs of generalized polygons show that
this is asymptotically tight when $k \in \{2,3,5\}$.
\end{abstract}

\section{Introduction}
In this note, we study graphs without cycles of prescribed even
lengths. For a finite or infinite set ${\cal C}$ of cycles, define
$\mbox{ex}(n,{\cal C})$ to be the maximum possible number of edges
in an $n$-vertex graph which does not contain any of the cycles in
${\cal C}$. The asymptotic behaviour of the function
$\mbox{ex}(n,{\cal C})$ is particularly interesting when at least
one of the cycles in ${\cal C}$ is of even length, and was
initiated by Erd\H{os} \cite{Erd}. In general, it is the lower
bounds for $\mbox{ex}(n,{\cal C})$ -- that is, the construction of
dense graphs without certain even cycles -- which are hard to come
by. The best known lower bounds are based on finite geometries,
such as polarity graphs of generalized polygons~\cite{LUW2}, and
the algebraic constructions given by Lazebnik, Ustimenko and
Woldar~\cite{LUW1} and Ramanujan graphs of Lubotsky, Phillips and
Sarnak~\cite{LPS}; see also~\cite{LUW3}. In the direction of upper
bounds, the first major result is known as the even circuit
theorem, due to Bondy and Simonovits \cite{BS}, who proved that
$\mbox{ex}(n,\{C_{2k}\}) \leq 100kn^{1+\frac{1}{k}}$. A more
extensive study of $\mbox{ex}(n,{\cal C})$ was carried out by
Erd\H{o}s and Simonovits~\cite{ES}. Our point of departure is the
study of $\mbox{ex}(n,{\cal C})$ when ${\cal C}$ consists only of
the even cycles of length at most $2k$. The main result of this
article is the following:

\begin{theorem}
\label{thm:main}
Let $k \geq 2$ be an integer. Then, for all $n$,
\[ \mbox{ex}(n,\{C_4,C_6,\dots,C_{2k}\}) \; \; \leq \; \;  \textstyle{\frac{1}{2}}n^{1 + \frac{1}{k}} + 2^{k^2}n.\]
Furthermore, when $k \in \{2,3,5\}$, the $n$-vertex
polarity graphs of generalized $(k + 1)$-gons in {\rm \cite{LUW2}}
have $\frac{1}{2}n^{1 + 1/k} + O(n)$
edges and no even cycles of length at most $2k$.
\end{theorem}

For the statement about the number of edges in the polarity
graphs, see~\cite{LUW2}, page 9. Theorem \ref{thm:main} extends
the Moore bound (see \cite{Big}) up to an additive term, and a
more recent result of Alon, Hoory, and Linial~\cite{AHL}, who
proved that an $n$-vertex graph without cycles of length at most
$2k$ has at most $\frac{1}{2}(n^{1 + 1/k} + n)$ edges (see
Proposition \ref{prop:AHL}). In other words, we do not require
that the odd cycles be forbidden, and the same bound still holds,
but with a weaker additive linear term. Our result is also best
possible in the following sense: if we forbid only the $2k$-cycle
in our graphs, then the upper bounds in Theorem \ref{thm:main} no
longer hold -- it was shown recently, in \cite{FNV}, that
$\mbox{ex}(n,\{C_6\}) > 0.534n^{4/3}$ and $\mbox{ex}(n,\{C_{10}\})
> 0.598n^{6/5}$ as $n$ tends to infinity.

\bigskip

\section{Local Structure}
Let $G$ be a graph with no even cycles of length less than or
equal to $2k$.  We write $P[u,v]$ to indicate that a path $P
\subset G$ has end vertices $u$ and $v$, and we order the vertices
of $P$ from $u$ to $v$. Let $\prec$ denote this ordering along $P$.
A {\it vine} on a path $P$ is a graph consisting of
the union of $P$ together with paths $Q[u_i,v_i]$
which are internally disjoint from $P$ for $i = 1,2,\dots,r$,
and where $u \preceq u_1 \prec v_1 \preceq u_2 \prec v_2 \preceq \dots
\preceq u_r \prec v_r \preceq v$.
A $uv$-path of shortest length is called a {\it $uv$-geodesic}.
A {\it $\theta$-graph} consists of three internally
disjoint paths with the same pair of endpoints.

\medskip

\begin{lemma}\label{theta}
Any $\theta$-graph contains an even cycle.
\end{lemma}

\begin{proof} If $P, Q$ and $R$ are the internally disjoint paths
in the $\theta$-graph with the same pair of endpoints, then
$|P \cup Q| + |Q \cup R| + |P \cup R| = 2|P| + 2|Q| + 2|R|$,
which is even. Therefore one of the cycles
$P \cup Q$, $Q \cup R$ or $P \cup R$ must have
even length.
\end{proof}

\medskip

\begin{lemma}
\label{lem:short}
Let $P^*$ be a $uv$-geodesic of length at most $k$. Then the union $H$
of all $uv$-paths of length at most $k$ is a vine on $P^{*}$ and $P^{*}$ is the unique $uv$-geodesic.
\end{lemma}

\begin{proof}
Suppose, for a contradiction, that $H$ is not a vine on $P^*$.
Let $x \prec v$ be a vertex of $P^*$ at a maximum distance from $u$ on $P^*$
such that the union of all $ux$-paths in $H$
is a vine on $P^*[u,x]$. By the maximality of $x$, there is a
$uv$-path $P$ of length at most $k$ such that $x$ has degree three in
$P \cup P^*$. If $P$ has minimum possible length, then
$P[x,y] \cup P^*[x,y]$ is the only cycle in $P \cup P^*$
for some $y \succ x$ on $P^*$. By the maximality of $x$, the union of all $uy$-paths
in $H$ is not a vine. Therefore there must be a $uv$-path $Q$ of
length at most $k$ such that $Q \cup P \cup P^*$ is not a vine on $P^*$.
If $Q$ has minimum possible length, then $P \cup Q$ and $P^* \cup Q$
each have exactly one cycle. It follows that there is a path $Q[w,z] \subset Q$ such that
\[ Q[u,x] = P^*[u,x] \; \;  \mbox{ and } \; \;  Q[x,w] \cup Q[z,v] \subset P[x,v] \cup P^*[x,v]\]

and $Q[w,z]$ is internally disjoint from $P \cup P^*$. Since $P
\cup P^* \cup Q$ is not a vine, $w \in P[x,y] \cup P^*[x,y]$ and
$w \neq y$. If $z \in P^*[y,v]$, then $P^*[x,z] \cup P[x,z] \cup
Q[w,z]$ is a $\theta$-graph (see Figure 1).

\medskip

\SetLabels
\R(.41*.00)$P^*[x,y]$\\
\R(.20*.44)$P^*$\\
\R(.35*.29)$x$\\
\R(.01*.29)$u$\\
\R(1.00*.29)$v$\\
\R(.82*.29)$z$\\
\R(.73*.29)$y$\\
\R(.90*.92)$Q[w,z]$\\
\R(.45*.65)$w$\\
\R(.65*.60)$P[x,y]$\\
\endSetLabels
\begin{center}
\centerline{\AffixLabels{\includegraphics[width=3in]{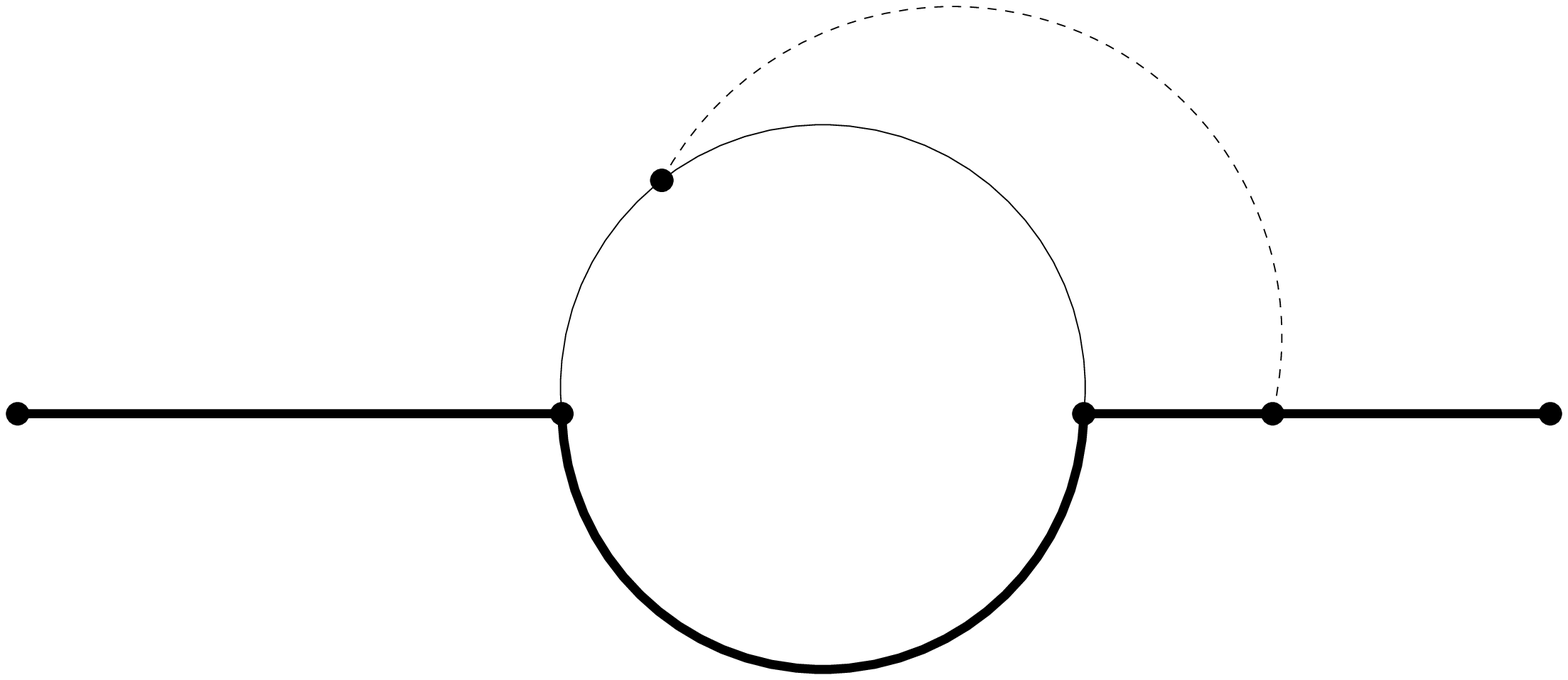}}}
\end{center}

\begin{center}
{\sf Figure 1 : A $\theta$-graph in $Q \cup P \cup P^*$.}
\end{center}

\medskip

The cycles in this $\theta$ graph are $P[w,z] \cup Q[w,z] \subset P \cup Q$
and $P[x,y] \cup P^*[x,y] \subset P \cup P^*$ and
$P^*[x,z] \cup Q[x,z] \subset P^* \cup Q$. Each
of these cycles has length at most $2k$, since the paths
$P,Q$ and $P^*$ each have length at most $k$.
By Lemma \ref{theta}, one of these cycles has even length, which is a contradiction.
A similar argument works when $z \not \in P^*[y,v]$.
Therefore $H$ is a vine on $P^*$.

\bigskip

To complete the proof, we must show that $P^{*}$ is the unique
$uv$-geodesic. By definition, $H$ consists of the union of $P^{*}$
and paths $P_i = P_i[u_i,v_i]$ for $i \in [r]$, and let $P^{*}_i =
P^{*}[u_i,v_i]$. Since each cycle $P^{*}_i \cup P_i$ is of length
at most $2k$, each cycle in the vine has odd length. Now suppose
$P$ is another $uv$-geodesic. Then $P_i \subset P$ for some $i$.
Since $P_i \cup P_i^{*}$ is an odd cycle, we may assume $|P_i| <
|P^{*}_i|$. By replacing $P^{*}_i$ with $P_i$ on $P^{*}$, we
obtain a $uv$-path of length $|P^{*}| - |P^*_i| + |P_i| <
|P^{*}|$, which contradicts the fact that $P^{*}$ is a
$uv$-geodesic. So $P^{*}$ is the unique $uv$-geodesic.
\end{proof}

\bigskip

Henceforth, the paths
in the vine on $P^{*}$ will be denoted $P_i = P_i[u_i,v_i]$,
and $P^{*}[u_i,v_i] = P^*_i$, for $i \in [r]$. Let $\p_k(u,v)$ denote
the set of all $uv$-paths of length $k$, and define the map
\[
f: \p_k(u,v) \rightarrow 2^{[r]} \; \; \mbox{ by } \; \;
f(P) = \set{i \in [r] \; \mid \; P_{i}[u_i,v_i] \subset P }.
\]
Then $f(P)$ records the set of integers $i$ for which the path $P
\in \p_k(u,v)$ uses the path $P_i[u_i,v_i]$ in the vine on $P^{*}$
instead of $P^{*}[u_i,v_i]$. Let ${\cal F}$ be the image of
$\p_k(u,v)$ under $f$.

\medskip

\begin{lemma}
\label{lem:paths}
The map $f$ is an injection, and the family ${\cal F}$
is an antichain of sets of size at most $k - |P^*|$ in
the partially ordered set of all subsets of $[r]$.
\end{lemma}
\begin{proof}
By Lemma \ref{lem:short}, each $P \in \p_k(u,v)$ is the union of
some (possibly none) of the paths $P_{i}$ together with internally
disjoint subpaths of $P^{*}$. Therefore the set $f(P)$ uniquely
determines $P$, and $f$ is an injection. If two sets in ${\cal F}$
are comparable, say $f(P) \subset f(Q)$, then $|Q| > |P|$ and $Q
\not \in \p_k(u,v)$, which is a contradiction. So ${\cal F}$ is an
antichain. Finally, any path $P \in \p_k(u,v)$ has length at least
$|P^{*}| + |f(P)|$, by Lemma \ref{lem:short}, so all sets in
${\cal F}$ have size at most $k - |P^{*}|$.
\end{proof}

\begin{theorem}\label{thm:max}
Let $G$ be a graph containing no even cycles of length at most $2k$.
Then
\[ |\p_k(u,v)| \; \leq \; \max\left( {r \choose m} : r \leq k \; \mbox{and}\; m = \min
\set{\Big\lfloor \frac{r}{2} \Big\rfloor,k-r}\right).\] The equality is achieved
when $r = |P^{*}|$ and the vine on $P^{*}$ comprises $|P^{*}|$
triangles.
\end{theorem}

\begin{proof}
The family ${\cal F}$ is an antichain, by Lemma \ref{lem:paths}.
By Sperner's Theorem and the LYM inequality \cite{Eng}, this means
that $|{\cal F}| \leq {r \choose m}$ where $m = \min
\set{\floor{\frac{r}{2}},k-|P^*|}$.
\end{proof}

\bigskip

A \emph{non-returning} walk of length $r$ in $G$ is a walk whose
consecutive edges are distinct. Let $\W_r$ be the set of
non-returning $r$-walks (for $r = 0$, $\W_0$ consists of single
vertices). The final result required for the proof of Theorem
\ref{thm:main} is the following lower bound on the number of
non-returning walks, by Alon, Hoory and Linial~\cite{AHL}, which
gives the best known upper bound on
$\mbox{ex}(n,\{C_3,C_4,\dots,C_{2k}\})$:

\begin{proposition}\label{prop:AHL}
Let $G$ be an $n$-vertex graph of average degree $d \geq 2$. Then
$|\W_r| \geq \; nd(d-1)^{r-1}$. Moreover, if $G$ has average degree $d \geq 2$ and no cycles of length at most
$2k$, then $d(d - 1)^{k-1} \leq n$.
\end{proposition}

In \cite{AHL}, the number $\W_r/nd$ is denoted $N_{r-1}$ and shown
to be less than $(d-1)^{r-1}$.  The second statement of the
Proposition is an immediate consequence of the main theorem there.
\bigskip

\section{Proof of Theorem \ref{thm:main}}

Let $G$ be a counterexample to Theorem \ref{thm:main} with minimal
number of vertices $n$ and average degree $d$. Then $d >
n^{\frac{1}{k}} + 2^{k^2}$, and $G$ has minimum degree at least
$\lfloor d/2 \rfloor + 1$, otherwise we remove a vertex of lower
degree, keeping the average degree non-increasing, to obtain a
smaller counterexample than $G$. We may also assume $n > 2^{k^2}$.
Now let $v$ be a vertex of $G$ of maximum degree, $\Delta$. Pick a breadth-first search tree $T$
rooted at $v$, and let $T_r$ be the set of vertices of $G$ at
distance at most $r$ from $v$. Then no vertex of $T_r$ is joined
to two vertices in $T_{r-1}$, and the set of edges in $T_{r-1}
\backslash T_{r-2}$ form a matching, for all $r \leq k$. So every
vertex of $T$ has degree at least $\delta - 2$, where $\delta$ is
the minimum degree in $G$, from which we deduce
\[ 1+\Delta +\Delta(\delta-2)+ \dots + \Delta(\delta -2)^{k-1} \; \leq \; |V(T)| \; \leq \; n.\]
Since $\delta > \lfloor d/2 \rfloor$ and
$d > n^{\frac{1}{k}} + 4$, we find $\Delta < 2^{k-1}n^{\frac{1}{k}}$.

\bigskip

Now let $\p_r$ be the set of paths of length $r$ in $G$, and let
$\q_r = \W_r - \p_r$ be the set of non-returning walks with $r$
edges which are not paths. There are at least $\delta - k$
extensions of a given path of length $r$ in $G$, for any $r < k$.
Therefore
\begin{equation}\label{eq:walkbound}
|\p_k| \geq  (\delta - k)^{k - \ell}|\p_{\ell}| \; \; \mbox{ and }
\; \; |\q_k| \leq \Delta^{k-1} k n < k2^{(k-1)^2}
n^{\frac{2k-1}{k}}.
\end{equation}
By Lemma \ref{lem:short}, for any pair $(u,v)$ of distinct
vertices, joined by at least two paths of length $k$, there is a
$uv$-geodesic of length $\ell < k$. By Theorem \ref{thm:max},
$|\p_k(u,v)| < 2^k$, so the number of ordered pairs of vertices
joined by exactly one $k$-path is at least
\begin{eqnarray*}
|\p_k| - 2^{k} \sum_{\ell = 1}^{k-1}|\p_{\ell}| &\geq& |\p_k|
\brac{ 1 - \frac{2^{k}}{\delta-k-1}}\\
&=& \left(\; |\W_k| - |\q_k| \;\right) \cdot \brac{ 1 - \frac{2^{k}}{\delta-k-1}}\\
&>& \left(nd(d-1)^{k-1} - k2^{(k-1)^2} n^{\frac{2k-1}{k}}\right)
\cdot \brac{ 1 - \frac{2^{k}}{\delta-k-1}}.
\end{eqnarray*}

In the last line, we used (\ref{eq:walkbound}) and Proposition
\ref{prop:AHL}. There are $n(n - 1)$ (ordered) pairs of distinct
vertices which could be joined by a unique path of length $k$, so
the expression above is less than $n^2$.  Using $\delta-k-1 \geq
\frac{d}{4}$ and substituting $d = n^{\frac{1}{k}} + 2^{k^2}$
into the last line, we get
\begin{eqnarray*}
n^2 &>& \left(n(n^{\frac{1}{k}} + 2^{k^2})(n^{\frac{1}{k}} +
2^{k^2} - 1)^{k-1} - k2^{(k-1)^2} n^{\frac{2k-1}{k}}\right)
\left(1 - \frac{2^{k+2}}{n^{\frac{1}{k}} + 2^{k^2}}\right) \\ \\
&=& \left(n^{\frac{2k - 1}{k}}(n^{\frac{1}{k}} + 2^{k^2})(1 +
n^{-\frac{1}{k}}(2^{k^2} - 1))^{k-1} -
k2^{(k - 1)^2}n^{\frac{2k - 1}{k}}\right)\left(1 - \frac{2^{k+2}}{n^{\frac{1}{k}} + 2^{k^2}}\right) \\ \\
&>& \left(n^{\frac{2k - 1}{k}}(n^{\frac{1}{k}} + 2^{k^2})(1 +
n^{-\frac{1}{k}}(k - 1)(2^{k^2} - 1)) -
k2^{(k - 1)^2}n^{\frac{2k-1}{k}}\right)\left(1 - \frac{2^{k+2}}{n^{\frac{1}{k}} + 2^{k^2}}\right) \\ \\
&>& n^2 \left(1 + \frac{2^{k^2}}{n^{\frac{1}{k}} + 2^{k^2}}\right)
\left(1 - \frac{2^{k + 2}}{n^{\frac{1}{k}} + 2^{k^2}}\right) \; \; > \; \;  n^2\\
\end{eqnarray*}
which gives a contradiction.  We must thus have $d <
n^{\frac{1}{k}} + 2^{k^2}$. \vrule height10pt width5pt depth1pt

\bigskip

\section{Concluding Remarks}

If $G$ is $d$-regular, then picking a breadth first search tree
as in the calculation of the maximum degree we obtain
         \[ 1+d+d(d-2)+ \dots +d(d-2)^{k-1} \leq n.\]
So in this case we have $d < n^{\frac{1}{k}} + 2$.
The main points at which the large linear term
is introduced in the proof of Theorem \ref{thm:main}
is in the estimate of the maximum degree
and the upper bound on $|{\cal Q}_k|$. We believe it
should be possible to circumvent these bounds to obtain
a linear term of the form $cn$, for some absolute constant $c$.
Finally, we note that the analogous extremal problem when some of
the short odd cycles are forbidden seems to be very difficult.
For example, it is known that
\[ \frac{1}{2\sqrt{2}} \; \leq \; \liminf_{n \rightarrow \infty} \frac{\mbox{ex}(n,\{C_3,C_4\})}{n^{3/2}} \; \leq \;
\limsup_{n \rightarrow \infty}
\frac{\mbox{ex}(n,\{C_3,C_4\})}{n^{3/2}} \; \leq\;  \frac{1}{2},\]
but the asymptotic value of $\mbox{ex}(n,\{C_3,C_4\})$ remains an
open question (posed by Erd\H{o}s).

\bigskip

\textbf{Acknowledgements.}
The first author would like to thank Terence Tao for supervising him during
his undergraduate thesis, which led to this work.

\end{document}